\documentclass[twocolumn]{autart}    

\usepackage[dvipsnames]{xcolor}
\usepackage{amsmath,amssymb}
\usepackage{mathrsfs}
\usepackage{algorithm}
\usepackage{algpseudocode}
\usepackage{float}
\usepackage[noadjust]{cite}
\usepackage{graphics} 
\usepackage{epsfig} 
\usepackage{float}
\usepackage{epstopdf}
\usepackage{graphicx}

\newtheorem{assumption}{Assumption}
\newtheorem{theorem}{Theorem}
\newtheorem{lemma}{Lemma}
\newtheorem{remark}{Remark}

\newtheorem{corollary}{Corollary}

\algrenewcommand\algorithmicrequire{\textbf{Input:}}
\algrenewcommand\algorithmicensure{\textbf{Output:}}

\begin{document}

\begin{frontmatter}

\title{Suboptimal Shrinking Horizon MPC with a Lower Hessian Condition Number from Adjustable Terminal Cost \thanksref{footnoteinfo}} 

\thanks[footnoteinfo]{Corresponding author Steven van Leeuwen}

\author[UMich,APL]{Steven van Leeuwen}\ead{svanlee@umich.edu},    
\author[UMich]{Ilya Kolmanovsky}\ead{ilya@umich.edu},               

\address[UMich]{Department of Aerospace Engineering, University of Michigan, Michigan, USA}  
\address[APL]{Johns Hopkins University Applied Physics Lab, Maryland, USA}             

\begin{keyword} Model Predictive Control, Shrinking Horizon, Computations             
\end{keyword}                             

\begin{abstract}                          
A strategy for reducing the number of iterations and computational burden in shrinking horizon Model Predictive Control (SH-MPC) when steering into a prescribed terminal set despite unmeasured disturbances is proposed. This strategy exploits dynamic adjustment of the terminal cost weight and horizon length while ensuring that the terminal set is reached within a desired number of steps. A lower Hessian condition number which facilitates the computational reduction is proved under assumptions, and an example of spacecraft nutation damping using the proposed approach is reported.
\end{abstract}

\end{frontmatter}

\section{Introduction}

Shrinking horizon Model Predictive Control (SH-MPC) is a variant of MPC where the Optimal Control Problem (OCP) is posed over a shrinking horizon instead of a fixed receding horizon; this can be advantageous in settings where terminal constraints need to be satisfied after a fixed amount of time, or in the setting of closed-loop trajectory optimization. In the SH-MPC setting, the impact of inexact computations, warm-starting, and unmeasured disturbances has been considered in \cite{SHMPC_main}. 

More specifically, for the projected gradient method (PGM), in \cite{SHMPC_main} bounds on the required number of solver iterations, which ensure that the trajectory is steered into the prescribed terminal set within a fixed amount of time, can be expressed as a function of the condition number of the Hessian. The condition number can be bounded in certain scenarios, such as when using prestabilization and when analyzing a block-Toeplitz infinite horizon matrix \cite{McInerney}. 

This paper proposes a strategy to reduce the condition number of the Hessian and in turn reduce the number of solver iterations through dynamically adjusting the terminal penalty weight and horizon length of the SH-MPC OCP. We also obtain new characterizations of the condition number behavior of the Hessian without preconditioning methods.

Notation: Given $a,b \geq 0$, let $\mathbb{Z}_{[a,b]} = \mathbb{Z} \cap [a,b]$. Given a symmetric matrix $W \succ 0$, the $W$-norm of x is $||x||_W = \sqrt{x^TWx}$. For any matrix $B$, we denote $||B||_W$ the max eigenvalue of $B^TWB$. Given $M \succ 0$, let $\lambda_{min}^W(M),\lambda_{max}^W(M)$ denote the min and max eigenvalues of $W^{-1/2}MW^{-1/2}$. Let the eigenvalues of a square matrix of size $N$ be ordered as $\lambda_{max} = \lambda_1 \geq \lambda_2 \geq ... \geq\lambda_{min}=\lambda_N$. Let $\textnormal{spec}$ denote the (ordered) spectrum of a square matrix of size $N$, i.e. $(\lambda_1,...,\lambda_N)$ and $\kappa(W) = \lambda_{max}(W)/\lambda_{min}(W)$.

\section{Preliminaries and Optimal Shrinking Horizon MPC}
\label{MPC and Preliminaries}
We consider a discrete-time Linear Time Invariant (LTI) system,
\begin{equation}\label{eq:dynamic_system}
    x_{k+1} = Ax_k + Bu_k + d_k, \ k \in \mathbb{Z}_{[0,N-1]},
\end{equation}
where $k$ is the time instant, $x_k \in \mathbb{R}^n$ is the state, $u_k \in \mathbb{R}^m$ is the control input, and $d_k \in \mathbb{R}^n$ is an unknown set-bounded disturbance, with $d_k \in \mathcal{D}$. The control objective is to reach a target terminal set $\Omega$ at the $N$th time instant so that $x_N \in \Omega$. There are also control constraints $ u_k \in \mathcal{U}$, where $\mathcal{U}$ is polyhedral. 

In SH-MPC, the control $u_k$ is generated by solving an OCP that traditionally shrinks so that its prediction horizon is of length $N-k$ steps at time instant $k$. In our SH-MPC formulation, the following OCP is solved at each time instant $k$, with the resulting $u_{k}$ applied to (\ref{eq:dynamic_system}):
\begin{subequations}
\begin{gather}\label{eq:OCP}
\mathcal{P}_{N-k}(x_k): \min_{\xi,\mu} \ \omega_k F(\xi_{N_k|k})+ \sum_{i=k}^{N_k-1} l(\xi_{i|k},\mu_{i|k})
\end{gather}
$s.t.$
\begin{gather}
\xi_{i+1|k} = A\xi_{i|k} + B\mu_{i|k}, \ i \in \mathbb{Z}_{[k,N_k-1]} \label{eq:OCP_b}\\
\mu_{i|k} \in \mathcal{U}, \ i \in \mathbb{Z}_{[k,N_k-1]} \label{eq:OCP_c}
\end{gather}
\end{subequations}
where $N_k \in \mathbb{Z}_{[N,\infty)},$ and $ \ \omega_k > 0$. Note that we do not explicitly impose the terminal constraint $\xi_{N|k} \in \Omega \label{eq:OCP_d}$
as a part of the OCP; however, this property will be ensured implicitly by later construction.
We define $\xi_{k|k} = x_k, \ z_k = (\mu_{k|k},...,\mu_{N_k-1|k})$,
\begin{equation}\label{eq:u_definition}
    u_k = \Xi_k z_k = \Xi_k(z^*_k(x_k) + e_k)  
\end{equation}
where $z^*_k(x_k)$ represents the optimal sequence corresponding to the solution of $\mathcal{P}_{N-k}(x_k)$, $z_k$ is an approximation of $z^*_k(x_k)$ with an error $e_k$ due to inexact computations, and $\Xi_k \in \mathbb{R}^{m \times (N_k-k)m}$ is a matrix that selects $\mu_{k|k}$ from $z_k$. We define $u_k^*(x_k) = \Xi_k z_k^*(x_k)$, and as in \cite{SHMPC_main}, focus on the case when $F(x) = ||x||_P^2, \ \Omega = \{ x \ | \ F(x) \leq \alpha\}, \ l(x,u) = ||x||_Q^2 + ||u||_R^2$. 

The following assumptions are made regarding the MPC formulation:
\begin{assumption}
The set $\mathcal{U}$ is closed, convex, and contains the origin in its interior.
\end{assumption}
\begin{assumption}
The pair $(A,B)$ is stabilizable.
\end{assumption}
\begin{assumption}
The weight matrices of (\ref{eq:OCP}) satisfy $Q \succ 0, R \succ 0$, and the Riccati equation, $P = Q + A^T P A -(A^TPB)K$, with $K = (R + B^TPB)^{-1} (B^TPA)$.
\end{assumption}

A tightened terminal set, dependent on a parameter $0 < \rho < 1$, is defined as $\Omega_{\rho} = \{ x \ | \ F(x) \leq \rho \alpha\}$, and the interior of $\Omega$ as $\Tilde{\Omega}$. We also define the $N_k-k$ backward reachable set of these terminal sets for system (\ref{eq:OCP_b}), (\ref{eq:OCP_c}) as $\Psi_{N-k}(\Omega), \ \Psi_{N-k}(\Omega)$ respectively. When $N_k = N$, we denote the backward reachable set as $\tilde{\Psi}_{N-k}(\Omega)$. The OCP $\mathcal{P}_{N-k}(x)$ can be written in the following condensed form (eliminating the state sequence) \cite{Liao-McPherson2022}
\begin{equation}\label{eq:OCP_condensed}
    \min_{z_k \in \mathcal{Z}_k} J_{N-k}(x,z_k) = \min_{z_k \in \mathcal{Z}_k}||(x,z_k)||^2_{M_k}
\end{equation}
where 
$\mathcal{Z}_k = \mathcal{U}^{N_k-k}$, and
$M_k = \begin{bmatrix}
    W_k & G_k^\top \\
    G_k & H_k
\end{bmatrix}$, see definitions of $W_k, \ G_k, \ H_k$ in \cite{Liao-McPherson2022}. We note $W_k, \ H_k \succ 0$. We also define $V_{N-k}(x) =  J_{N-k}(x,z^*_k(x)), \ \psi_{N-k}(x) = \sqrt{V_{N-k}(x)}$.

\subsection{Regularity Properties}
\begin{lemma} (\cite{Liao-McPherson2022} Corollary 2)
    Let Assumptions 1-3 hold, then $\forall k \in \mathbb{Z}_{[0,N-1]}$ and $\forall x, \ y \in \mathbb{R}^n$, the solution mapping $x \mapsto z^*_k(x)$ satisfies
    \begin{gather*}
        \langle z_k^*(x) - z_k^*(y),G_k(x-y) \rangle \leq -||z_k^*(x)-z_k^*(y)||_{H_k}^2, \\
        ||z_k^*(x)-z_k^*(y)||_{H_k} \leq ||G_k(x-y)||_{H_k^{-1}} \\
        |\psi_{N-k}(x)-\psi_{N-k}(y)| \leq ||x-y||_{W_k}.
    \end{gather*} 
    From \cite{Liao-McPherson2022}, Corollary 2, $z^*_k(x)$ is Lipschitz continuous.
\end{lemma}

\section{Dynamic Adjustment of the OCP}
In our proposed approach, the horizon length $N_k$ and terminal penalty weight $\omega_k$ are dynamically adjusted to allow for better conditioning of the Hessian $H_k$. We prove in Section 4.3 that, under suitable assumptions, the condition number of $H_k$ decreases with a lower value of $\omega_k$ for an unchanged horizon length. Increasing the horizon length and lowering the value of $\omega_k$ can also decrease the condition number of $H_k$, but further analysis of this is left to future work. When $\omega_k$ is adjusted, \cite{stable_no_terminal_constraint} shows a backward reachable set estimate (BRSE) of $\Omega$ is altered; this estimate informs the choice of $N_k, \ \omega_k$.

We have the following assumption to ensure positive invariance of the terminal set $\Omega$:
\begin{assumption}(\cite{SHMPC_main} Assumption 7)
   \begin{equation} \label{eq:invariance_property}
    \min_{u \in \mathcal{U}}\{F(Ax+Bu) - F(x) + l(x,u)\} \leq 0, \ \forall x \in \Omega.
\end{equation} 
\end{assumption}

With this, a BRSE when solving $\mathcal{P}_{N-k}(x_k)$ is  $\Gamma_{N-k} = \{ x \ | \ V_{N-k}(x) \leq l(x,u_k^*(x)) + (N-k-1)\gamma + \omega_k \alpha \}$, with $\gamma = \alpha \lambda_{min}^P(Q)$, which ensures that \cite{SHMPC_main} $\xi_{i|k} \in \Omega, \ i \in \mathbb{Z}_{[N,N_k]}$. A tighter BRSE is  
\begin{equation}\label{eq:BRSE_estimate}
    \Upsilon_{N-k} = \{ x \ | \ V_{N-k}(x) \leq (N-k)\gamma + \omega_k \alpha \} \subseteq \Gamma_{N-k}.
\end{equation}

When the dependence on $\omega_k$ is emphasized in (\ref{eq:BRSE_estimate}), we write $\Upsilon_{N-k}(\omega_k).$ Then the following theorem underpins our methodology for constructing $\omega_k$.

\begin{theorem} {(\cite{stable_no_terminal_constraint} Theorem 3)}
    Consider $\mathcal{P}_{N-k}(x)$. For any $x \in \tilde{\Psi}_{N-k}(\tilde{\Omega})$ there exists a finite $\omega_k > 0$ such that $x \in \Upsilon_{N-k}(\omega_k)$.
\end{theorem}

Theorem 1 is satisfied by letting $\omega_k \geq \Lambda_k(z_{k},\rho_{k})$ \cite{stable_no_terminal_constraint}, where
\begin{equation}\label{eq:Lambda}
    \Lambda_k(z_{j},\rho) = \text{max}\bigg(\epsilon, \ \frac{L_{N-k}(z_j) - (N-k)\gamma}{(1-\rho)\alpha}\bigg),
\end{equation} 
and where $\epsilon>0$ is a small number close to zero, $L_{N-k}(z_j) = \sum_{i \in [k,N-1]} l(\mu_{i|j},\xi_i(z_j)), \ j \leq k$, $\xi_k(z_j) = x_k, \
 \xi_{i+1}(z_j) = \xi_i(z_j) + B\mu_{i|j}, \ i \in \mathbb{Z}_{[k,N-2]}$, (a possibly suboptimal sequence), and $\xi_N(z_j) \in \tilde{\Omega}$.
\begin{remark}
    For $z_j, \ \rho$ defined in (\ref{eq:Lambda}), from (\ref{eq:invariance_property}) a sequence with a longer horizon is feasible, which if used in (\ref{eq:Lambda}), leads to a reduced computed value of $\rho$, which in turn can reduce the terminal penalty weight.
\end{remark}

\section{Suboptimality}
\subsection{Value Function Decay}
 We establish conditions such that $x_k \in \Upsilon_{N-k}$ under the dynamically adjusted OCP, inexact computations, and a disturbance $d_k$. This property ensures $x_N \in \Omega$ and $\rho_k < 1$.
 
 We first describe a difference between the cost function value and the sublevel set value in $\Gamma_{N-k}$, defined as $\beta_k'$, which we compute when the parameters of (\ref{eq:OCP}) change, i.e., $N_k, \ \omega_k$. Let $\Delta k$ be the difference between $k > 0$ and the last time instant $N_k$ and/or $\omega_k$ were changed.
\begin{assumption}
    If $\Delta k = 0$, then for $\omega_k = \Lambda_k(z_{k-1},\rho_{k})$, there exists $\beta'_k >0$ such that $V_{N-k}(x_k) \leq l(x_k,u^*_k) + (N-k-1)\gamma + \omega_k \alpha - \beta_k'$.
\end{assumption}
Assumption 5 provides a margin for the additional cost accrued by suboptimal cost functions $J_{N-k}(x_k,z_k)$, while still ensuring $x_k \in \Upsilon_{N-k}$. Then, considering the exact system and inexact system,
\begin{gather}
    x_{k+1}^* = Ax_k + B u^*_k(x_k), \ x^*_0 = x_0, \label{eq:xp1_star}\\
    x_{k+1} = Ax_k + Bu^*_k(x_k) + w_k, \label{eq:xp1}
\end{gather}
where from (\ref{eq:dynamic_system}), (\ref{eq:u_definition}), $w_k = d_k + B \Xi_ke_k$, we have the following:
\begin{lemma}
Let Assumptions 1-3 hold, and $\omega_{k+1} \leq \omega_k, \ N_{k+1} \geq N_k$. Suppose $x_{k} \in \tilde{\Psi}_{N-k}(\Omega)$, and $x_{k+1}$ is given by (\ref{eq:xp1}), then
\begin{gather}
    \psi_{N-k-1}(x_{k+1}) \leq \psi_{N-k-1}(x^*_{k+1}) + ||w_k||_{W_{k+1}}, \label{eq:Lemma6_1}\\
    \psi_{N-k-1}(x^*_{k+1}) \leq \sqrt{V_{N-k}(x_k) - l(x_k,u_k^*(x_k))} . \label{eq:Lemma6_2}
\end{gather}
\end{lemma}
\begin{pf}
    The proof follows that in (\cite{SHMPC_main} Proposition 14). Equation (\ref{eq:Lemma6_1}) follows from Lemma 1. To prove (\ref{eq:Lemma6_2}), $\psi_{N-k-1}(x^*_{k+1}) \leq \{ \psi_{N-k-1}(x^*_{k+1}) \ | \ \omega_{k+1} = \omega_k, \ N_{k+1} = N_{k} \}$ which in turn is equal to the right hand side in (\ref{eq:Lemma6_2}) by the principle of optimality. \qed
\end{pf}

We next present the following result which formalizes satisfying the terminal constraint with the update (\ref{eq:xp1}) given $\omega_k$.
\begin{theorem} 
    Let $\Delta k = 0$, Assumptions 1-5 hold, and consider the closed loop system given by (\ref{eq:xp1}) at time instant $k$. Let $\{\beta_i\}_{i\in[k,N]}$ be any sequence satisfying $\beta_k \leq \beta_k'$, with $\beta_i \in (0,(N-i-1)\gamma + \omega_i \alpha], \forall i \in \mathbb{Z}_{[k,N-1]}, \ \beta_N = 0$, and let the following hold for all $i \in \mathbb{Z}_{[k,N-1]}$,
    \begin{gather*}
	\begin{array}{ll}
		(\omega_i - \omega_{i+1})\alpha < l(x_i,u_i^*(x_i)) - \gamma,  & \mbox{if} \ x_i \notin \Omega \\
		\omega_{i+1} = \omega_i,  & \mbox{if} \ x_i \in \Omega
	\end{array}
    \end{gather*}
    with $\omega_{i+1} \leq \omega_i, \ N_{i+1} \geq N_i, \ \omega_{N} = \omega_{N-1}$. Then define 
    \begin{gather*}
        \bar{V}_i = (N-i)\gamma + \omega_{i} \alpha - \beta_i \\
        \bar{w}_i = \pi_i^{-1}(\bar{V}_{i+1}^{1/2} - [\bar{V}_{i+1} - (\beta_i - \beta_{i+1})]^{1/2})
    \end{gather*}
    where $\pi_i = \lambda_{max}^{1/2}(W_{i+1})$ for $i \in \mathbb{Z}_{[k,N-2]}, \ \pi_{N-1} = \sqrt{\omega_{N-1} \lambda_{max}(P)}$. If $||w_i|| \leq \bar{w}_i$ for all $i \in \mathbb{Z}_{[k,N-1]}$, then 
    \begin{gather}
        V_{N-i}(x_i) \leq \bar{V}_i, \forall i \in \mathbb{Z}_{[1,...,N-1]}, \label{eq:Theorem2_1}\\
        F(x_N) \leq \alpha. \label{eq:Theorem2_2}
    \end{gather}
    Moreover, $x_i \in \Upsilon_{N-i}(\omega_i)$ for all $i \in \mathbb{Z}_{[k,...,N-1]}, u_i \in \mathcal{U}$ for all $i \in \mathbb{Z}_{[k,...,N-1]}$ and $x_N \in \Omega$.
\end{theorem}
\begin{pf}
    The proof follows that of (\cite{SHMPC_main}, Theorem 15). \qed
\end{pf}
The following establish bounds for which (\ref{eq:Theorem2_1}),(\ref{eq:Theorem2_2}) hold.
\begin{assumption}{(\cite{SHMPC_main} Assumption 18)}
There exists a known constant $\bar{d} > 0$ satisfying $\{d \ | \ ||d|| \leq \bar{d}\} \supseteq \mathcal{D}$ and $\bar{d} \leq \bar{w}_i \ \forall i \in \mathbb{Z}_{[k,N-1]}$.
\end{assumption}
\begin{corollary}
    Let Assumption 6 hold and define $\bar{e}_i = (1/||B||) (\bar{w}_i - \bar{d})$. If $||e_i|| \leq \bar{e}_i \ \forall i \in \mathbb{Z}_{[k,N-1]}$, then (\ref{eq:Theorem2_1}),(\ref{eq:Theorem2_2}) hold under the conditions in Theorem 2.
\end{corollary}

\subsection{Online Algorithm}
The following procedure is an example of how to dynamically adjust the OCP under Theorem 1 and 2. We let $N_k = N$, which guarantees a lower Hessian condition number at each time instant under Assumption 7 to follow. The online algorithm is driven by ease of implementation since $\{\bar{w}_k\}$ in Theorem 2 cannot be precomputed offline without also precomputing $\{\omega_k\}$.

i.) ($k=0$) We choose $\omega'_0$ such that is satisfies Assumption 5 and solve for $V_N(x_0)$. Let $\tilde{\omega}_0 \in [\Lambda_0(z^*_0(x_0),\rho_0),\ \omega'_0]$, otherwise $\omega_0 = \omega'_0$. Compute the resulting $\tilde{J}_{N}$ using $z^*_0(x_0)$ and $\tilde{\omega}_0$, see (\ref{eq:OCP_condensed}). Go to iii.)

ii.) ($k>0$) We check if $(\mu_{k|k-1},...,\mu_{N-1|k-1})$ is feasible. If so, compute the resulting $\tilde{\rho}_k$, otherwise $\omega_k = \omega_{k-1}$. Let $\tilde{\omega}_k \in [\Lambda_k(z_{k-1},\tilde{\rho}_k), \ \omega_{k-1}]$, otherwise $\omega_k = \omega_{k-1}$. Compute the resulting $\tilde{J}_{N-k}$ using $(\mu_{k|k-1},...,\mu_{N-1|k})$ and $\tilde{\omega}_k$. Go to iii.)

iii.) Compute $\beta_k' = \tilde{J}_{N-k} -(N-k)\gamma - \tilde{\omega}_k \alpha$. If $\beta'_k > 0$, then $\omega_k = \tilde{\omega}_k$ and $\{\bar{w}_i\}_{i \in [k,N]}, \ \{\bar{e}_i\}_{i \in [k,N]}, \ \bar{d}$ are recomputed subject to Theorem 2, assuming $\omega_{i+1} = \omega_{i}, \  \forall i \in \mathbb{Z}_{[k,N-1]}$, and Assumption 6. Otherwise, $\omega_k = \omega_{k-1}$ or for $k=0$, $\omega_0 = \omega'_0$.

\subsection{Optimizer and Iteration Bounds}

We have shown that the terminal penalty weight can be lowered per Theorem 1. Now we set out to show a reduction in computational burden that results from this as opposed to holding the terminal weight constant for $k \in \mathbb{Z}_{[0,N-1]}$. First, the quadratic form of the Hessian $H_k$ can be written as
\begin{equation}\label{eq:Hessian_expanded}
    ||v||_{H_k}^2 = \omega_k ||S_{N_k}v||_P^2 + ||v||_{\bar{R}}^2 + \sum_{l\in [1,N_k-1]} ||S_lv||_Q^2, 
\end{equation}
where $H_k = S^T\bar{Q}S + \bar{R}, \ S = (S_1,...,S_{N_k})$, $S_l = [A^{l-1}B, \ ..., \ AB, \ B, \ 0, \ ..., \ 0] \in \mathbb{R}^{n \times N_km}$, and $\bar{Q} = \textnormal{blkdiag}(I_{N_k-1} \otimes Q, \ \omega_k P), \ \bar{R} = \textnormal{blkdiag}(I_N \otimes R$).
    
We will need the following assumption for time instants $k=0$ through $N_{\kappa} < N$:
\begin{assumption}
The following is satisfied.
\begin{equation}\label{eq:Assumption_7}
    \frac{||B||_P}{||R||} \leq \frac{||S_{N_k}\overline{v}_k^0||_P^2}{||\overline{v}_k^0||_{H^1_k}^2}, \ \forall k \in \mathbb{Z}_{[0,N_{\kappa}]},
\end{equation}
and $\lambda_{max}(H_k), \lambda_{min}(H_k)$ are unique $ \forall k \in \mathbb{Z}_{[0,N_{\kappa}]}, \forall \omega \in (0,\omega'_0]$, and $R = \chi I$, where $\chi >0$ is a scalar,
\end{assumption}
where we define normalized eigenvectors 
\begin{gather*}
    \underline{v}_k^0 = \arg \min_{||v|| = 1} ||v||_{H^0_k}^2, 
    \overline{v}_k^0 = \arg \max_{||v|| = 1} ||v||_{H^0_k}^2, \\
    \underline{v}'_k = \arg \min_{||v|| = 1} ||v||_{H_k'}^2, \overline{v}'_k = \arg \max_{||v|| = 1} ||v||_{H_k'}^2,
\end{gather*}
where $H_k^1 = H_k(\omega_k = \omega'_0), \ H_k^0 = H_k(\omega_k = 0)$, and where $H'_k$ is associated with $\{ \mathcal{P}_{N-k}(x_k) \ | \ \omega_k = \omega_k' \}$. We note the uniqueness of the max and min eigenvalues in Assumption 7 is a technical assumption which can be relaxed, but this is not the focus of the paper. Now we state the main result.
\begin{theorem}
    Assume $x_k \in \Upsilon_{N-k}(\omega'_k)$. Then for $\omega_k < \omega'_k, \ \forall k \in \mathbb{Z}_{[0,N_{\kappa}]}$, 
    \begin{equation}\label{eq:spectral_reduction}
        H_k \prec H_k', \ \textnormal{spec}(H_k) < \textnormal{spec}(H_k'),
    \end{equation}
    where the second inequality is in an element-wise sense. Furthermore, if Assumption 7 holds,
    \begin{equation}\label{eq:condition_num_reduction}
        \kappa(H_k) < \kappa(H'_k). 
    \end{equation}
\end{theorem}
\begin{pf}
    We denote the optimal controls sequence and optimal cost associated with $\{ \mathcal{P}_{N-k}(x_k) \ | \ \omega_k = \omega_k' \}$ as $z_k'^*(x_k), \ V_{N-k}'(x_k)$, respectively. Noting each term in (\ref{eq:Hessian_expanded}) is positive definite, we have $H_k \prec H_k'$. 
    
    To show the spectral reduction in (\ref{eq:spectral_reduction}), we use the Courant-Fischer Minimax Theorem (\cite{Golub}, Page 394)
    \begin{equation}
        \lambda_l(H_k) = \max_{dim(\mathcal{V}) = l} \big(\min_{0 \neq v \in \mathcal{V}} \frac{|| v||_{H_k}^2}{||v||^2}\big), \ l \in \mathbb{Z}_{[1,...,N_km]},
    \end{equation}
    where $\mathcal{V}$ is any subspace of $\mathbb{R}^{N_km}$. Then $\lambda_l(H_k) < \lambda_l(H'_k)$ follows from $H_k \prec H_k'$.

    To show (\ref{eq:condition_num_reduction}), we first aim to show 
    \begin{equation}\label{eq:kappa_decrease_necessary}
        \frac{||S_{N_k}\underline{v}'_k||_P^2}{||\underline{v}'_k||_{H^1_k}^2} < \frac{||S_{N_k}\overline{v}'_k||_P^2}{||\overline{v}'_k||_{H^1_k}^2}, \ \forall k \in \mathbb{Z}_{[0,N_{\kappa}]}, \ \forall \omega_k' \in (0,\omega'_0],
    \end{equation}
    which we show through the following series of inequalities.
    \begin{gather}\label{eq:lhs_bound_kappa_rearrangement}
    \frac{||S_{N_k}\underline{v}'_k||_P^2}{||\underline{v}'_k||_{H^1_k}^2} < \frac{||B||_P}{||R||} \leq \frac{||S_{N_k}\overline{v}_k^0||_P^2}{||\overline{v}_k^0||_{H^1_k}^2} \leq \frac{||S_{N_k}\overline{v}'_k||_P^2}{||\overline{v}'_k||_{H^1_k}^2},
    \end{gather}
    where the first inequality is shown as follows. First, we have $||R|| = ||\underline{v}'_k||^2_{\bar{R}} < ||\underline{v}'_k||_{H^1_k}^2$. Secondly, note that the normalized vector $v$ from (\ref{eq:Hessian_expanded}) can be chosen such $||S_{N_k}v||_P^2 = ||B||_P$, then $||v||^2_{H'_k} = \omega'||B||_P + ||R||$. This, combined with $||\underline{v}'_k||^2_{H'_k} \leq ||v||^2_{H'_k}$ from (\ref{eq:Hessian_expanded}) imply $||S_{N_k}\underline{v}'_k||_P^2 \leq ||B||_P$. Combining these two observations the inequality follows. 
    
    The second inequality in (\ref{eq:lhs_bound_kappa_rearrangement}) is Assumption 7. 
    
    The third inequality in (\ref{eq:lhs_bound_kappa_rearrangement}) is shown as follows. We have $\Delta^2_Q = \sum_{l \in [1,N_k-1]}||S_{l}\overline{v}^0_k||_Q^2 - \sum_{l \in [1,N_k-1]}||S_{l}\overline{v}'_k||_Q^2 \geq 0$ (\ref{eq:Hessian_expanded}) since $ ||\overline{v}_k^0||_{H_k^0}^2 \geq ||\overline{v}_k'||_{H_k^0}^2$ by definition of $\overline{v}_k^0, \ \overline{v}_k'$ and by $\bar{R}$ being diagonal. Then $||S_{N_k}\overline{v}^0_k||_P^2 \leq ||S_{N_k}\overline{v}'_k||_P^2$ follows from $||\overline{v}_k^0||_{H_k'}^2 \leq ||\overline{v}'_k||_{H_k'}^2$. Letting $\Delta^2_P = \omega'_0||S_{N_k}\overline{v}'_k||_P^2 - \omega'_0||S_{N_k}\overline{v}^0_k||_P^2$, then
    \begin{gather}\label{eq:third_inequality_expanded}
         \frac{\omega'_0||S_{N_k}\overline{v}^0_k||_P^2}{||\overline{v}^0_k||_{H^1_k}^2 } \leq \frac{\omega'_0||S_{N_k}\overline{v}^0_k||_P^2 + \Delta^2_P}{||\overline{v}^0_k||_{H^1_k}^2 + \Delta^2_P - \Delta^2_Q} = \frac{\omega'_0||S_{N_k}\overline{v}'_k||_P^2}{||\overline{v}'_k||_{H^1_k}^2}. 
    \end{gather}
    By noting the first inequality in (\ref{eq:third_inequality_expanded}) holds since $\omega'_0||S_{N_k}\overline{v}^0_k||_P^2 \leq||\overline{v}^0_k||_{H^1_k}^2 $ and the denominator was increased by a smaller amount than the numerator, the third inequality in (\ref{eq:lhs_bound_kappa_rearrangement}) directly follows.

    Then from (\ref{eq:kappa_decrease_necessary}) we divide both sides by $||\overline{v}_k'||^2_{H'_k}/||\overline{v}_k'||^2_{H^1_k}$, and observing $||\overline{v}_k'||^2_{H^1_k} / ||\overline{v}_k'||^2_{H'_k} > 1$, if follows
        \begin{equation}\label{eq:kappa_decrease_necessary_full}
        \frac{||S_{N_k}\underline{v}_k'||_P^2}{||\underline{v}_k'||^2_{H'_k}} < \frac{||S_{N_k}\overline{v}_k'||_P^2}{||\overline{v}_k'||^2_{H'_k}}.
    \end{equation}
Observing that $\partial \lambda_{min}(H_k')/\partial{\omega_k'} = {\underline{v}'_k}^T (\partial(H'_k)/\partial \omega'_k) {\underline{v}'_k},$ and $ \partial \lambda_{max}(H_k')/\partial{\omega_k'} = {\overline{v}'_k}^T (\partial(H'_k)/\partial \omega'_k) {\overline{v}'_k}$, from (\ref{eq:Hessian_expanded}) we have $\partial \lambda_{min}(H_k')/\partial{\omega_k'} = ||S_{N_k}\underline{v}_k'||_P^2$ and $ \partial \lambda_{max}(H_k')/\partial{\omega_k'} = ||S_{N_k}\overline{v}_k'||_P^2$. Thus we rewrite (\ref{eq:kappa_decrease_necessary_full}) as
\begin{equation}\label{eq:chain_rule_sufficient_cond}
 \frac{\partial \lambda_{min}(H'_k)/\partial\omega_k'}{ \lambda_{min}(H'_k)} < \frac{\partial \lambda_{max}(H'_k)/ \partial\omega_k'}{\lambda_{max}(H_k')}.    
\end{equation}
We then prove (\ref{eq:condition_num_reduction}) by claiming 
\begin{equation} \label{eq:chain_rule_here}
    \partial \kappa(H'_k)/\partial \omega_k' > 0 \ \forall \omega'_k \in (0,\omega'_0].
\end{equation} After applying the chain rule to the l.h.s of (\ref{eq:chain_rule_here}) and simplifying, a sufficient condition for (\ref{eq:chain_rule_here}) to hold is (\ref{eq:chain_rule_sufficient_cond}). The function $\kappa(H'_k)$ is continuously differentiable on $(0,\omega'_0]$ from Assumption 7. \qed
\end{pf}

We next derive an offline bound on the number of iterations $\bar{l}_k$ of the PGM. See (\cite{SHMPC_main}, Section 6) for the derivation of (\ref{eq:l_bound}), and \cite{Liao-McPherson2022} for the PGM implementation. The offline bounds are, for all $k \in \mathbb{Z}_{[0,...,N-1]}$,
\begin{subequations} \label{eq:l_bound}
\begin{gather} 
    \begin{split}
    \bar{l}_{k} = \lceil\log \left( \frac{\lambda_{min}^{1/2}(H_k)\bar{e}_k}{||G_k x_k||_{H_k^{-1}}} \right) / \log \left( \frac{\kappa(H_k)-1}{\kappa(H_k)+1} \right) \rceil, \ \Delta k = 0 \label{eq:l_0_bound} \\
    \end{split} \\
    \begin{split} 
    \bar{l}_{k} = \lceil \log \left( \frac{\bar{e}_k}{(1+\tau_k)\bar{e}_{k-1} + \sigma_k \bar{d}} \right) / \log \left( \frac{\kappa(H_k)-1}{\kappa(H_k)+1} \right) \rceil, \\ \Delta k \neq 0 \label{eq:l_k_bound}
    \end{split}
\end{gather}
\end{subequations}
where $\tau_k = \lambda_{min}^{-1/2}(H_k)||{H_k}^{-1/2} G_k B||$ and where $\sigma_k = \lambda_{min}^{-1/2}(H_k)||{H_k}^{-1/2}G_k||$.

\section{Numerical Example}
\begin{figure}
\hspace*{-0.1cm}\vspace*{-0.0cm}\includegraphics[scale=0.34]{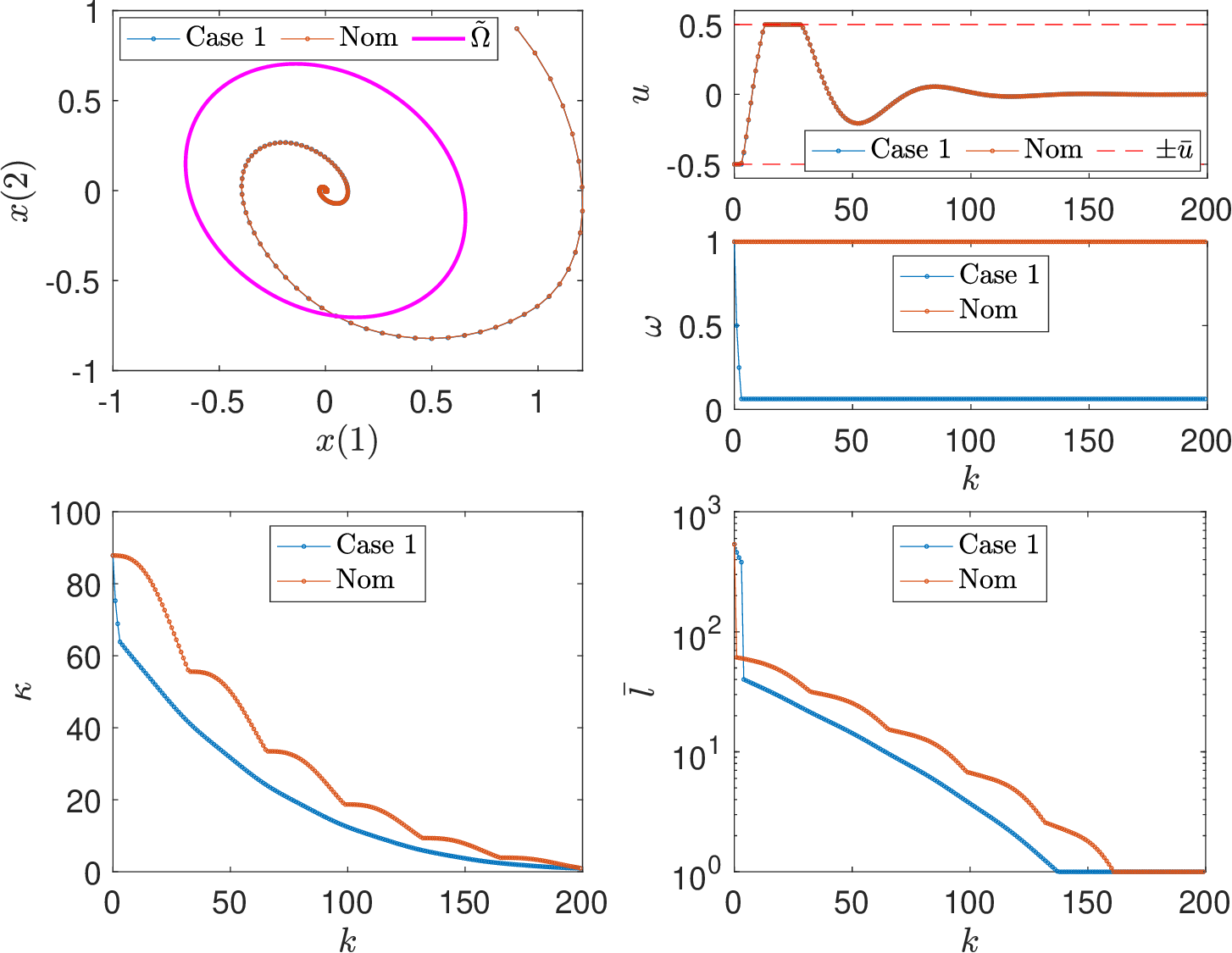}
\label{2}
\caption{State and control Trajectories for the axisymmetric spacecraft spin stabilization example.}
\end{figure}

We compare the online algorithm of this paper (called case 1) to the offline implementation in \cite{SHMPC_main} where $\omega_k = \omega_0'$ for all $k$ (called the nominal case), when using iteration bounds (\ref{eq:l_bound}) and confirming Assumption 7 holds. We use the axisymmetric spacecraft spin stabilization example from \cite{SHMPC_main}. The continuous-time equations of motion for the transversal components $q$ and $r$ of the angular velocity vector are
\begin{equation*}
    \dot{q} = ar + M_d/E_y, \ \dot{r} = -aq + M_c/E_y
\end{equation*}
where $p = 1, a = p(E_y - E_x)/E_y, \ E_y = 1, \ E_x = 0.05$. See \cite{SHMPC_main} for definitions of each of these terms. The disturbance input is $M_d$. We have $x = (q,r), u = M_c$ with $\mathcal{U} = [-0.5,0.5]$. The system is discretized using a sampling period of $0.1$, $N = 200$, $Q = \text{diag}(1,1), \ R= 3, \ \omega_0' = 1, \ \alpha = 13.2667$. An initial condition of $x_0 = (0.9,0.9)$ is used. The margin $\beta'_0$ for both cases is $\beta'_0 = l(x_0,u^*_0(x_0)) + (N-1)\gamma + \omega_0\alpha - V_N(x_0)$. For the nominal case, $\bar{d} = 1.58 \times 10^{-4}$, for case 1 $\bar{d}$ took values $[1.58,1.31,1.31,1.36] \times 10^{-4}$. In Figure 1, the blue graphs represent case 1 and orange graphs the nominal case. In the upper left plot, both cases show $x_N$ is in the terminal set. The condition number of the Hessian, seen in the lower left plot, is strictly lower for case 1, consistent with Proposition 1, since the $\omega_k$ plot shows values of the terminal penalty weight strictly lower for case 1. The iteration bounds are lower everywhere for case 1 except the first four time instants when the zero vector warmstart is used (\ref{eq:l_0_bound}), as seen in the lower right plot. To show a reduction in computational burden, we sum up the flop count for each iteration at each time instant. The expression for the flop count for one iteration at time instant $k$ is given by $(N-k)m(2((N-k)m)-1) + (N-k)m(2n-1) + 5(N-k)m$. The total flop count is still reduced (case 1: $1.24 \times 10^7$, nominal case: $2.19 \times 10^7$) in spite of the instances of zero warmstart.
\bibliographystyle{plain}        \bibliography{shmpc}

\end{document}